\documentclass[12pt]{article}

\usepackage{graphicx}
\usepackage{epsfig}

\usepackage[T1]{fontenc}
\input cyracc.def

\font\twlcyr=wncyr10 at 12pt

\def\cyr{\twlcyr\cyracc}

\usepackage{latexsym}
\newsavebox{\toy}
\savebox{\toy}{\framebox[0.65em]{\rule{0cm}{1ex}}}

\newcommand{\QED}{\usebox{\toy}\end{demo}}
\newenvironment{property}%
{\begin{list}{}{\setlength{\rightmargin}{0pt}%
\setlength{\itemsep}{0pt}}}{\end{list}}
\newlength{\templength}
\newcommand{\bp}{\setlength{\templength}{\labelwidth}%
\setlength{\labelwidth}{2em}\begin{property}}
\newcommand{\ep}{\end{property}\setlength{\labelwidth}{\templength}}
\newtheorem{theorem}{Theorem}[section]
\newtheorem{lemma}[theorem]{Lemma}
\newtheorem{proposition}[theorem]{Proposition}
\newtheorem{corollary}[theorem]{Corollary}
\newtheorem{assumption}{Assumption}
\newtheorem{definition}{Definition}[section]
\newtheorem{remark}{Remark}[section]
\newtheorem{exercise}{Exercise}[section]
\newtheorem{example}{Example}[section]

\newcommand{\Theorem}[1]{\begin{theorem}\label{Thm.#1}}
\newcommand{\Lemma}[1]{\begin{lemma}\label{Lem.#1}}
\newcommand{\Proposition}[1]{\begin{proposition}\label{Prop.#1}}
\newcommand{\Corollary}[1]{\begin{corollary}\label{Cor.#1}}
\newcommand{\Assumption}[1]{\begin{assumption}\label{Ass.#1}\rm}
\newcommand{\Definition}[1]{\begin{definition}\label{Def.#1}\rm}
\newcommand{\Remark}[1]{\begin{remark}\label{Rem.#1}\rm }
\newcommand{\Exercise}[1]{\begin{exercise}\label{Exe.#1}\rm }

\def\qed{\hfill\rule{.2cm}{.2cm}\par\medskip\par\relax}

\newcommand{\bd}{\begin{displaymath}}
\newcommand{\ed}{\end{displaymath}}
\newcommand{\bdn}{\begin{equation}}
\newcommand{\bdnl}{\begin{equation}\label}
\newcommand{\edn}{\end{equation}}
\newcommand{\barray}{\begin{array}}
\newcommand{\earray}{\end{array}}
\newcommand{\bds}{\begin{description}}
\newcommand{\eds}{\end{description}}
\newcommand{\bitemize}{\begin{itemize}}
\newcommand{\eitemize}{\end{itemize}}
\newcommand{\benumerate}{\begin{enumerate}}
\newcommand{\eenumerate}{\end{enumerate}}
\newcommand{\btabbing}{\begin{tabbing}}
\newcommand{\etabbing}{\end{tabbing}}
\newcommand{\bcenter}{\begin{center}}
\newcommand{\ecenter}{\end{center}}
\newcommand{\bflushright}{\begin{flushright}}
\newcommand{\bflushleft}{\begin{flushleft}}
\newcommand{\eflushright}{\end{flushright}}
\newcommand{\eflushleft}{\end{flushleft}}
\newcommand{\bdnn }{\begin{eqnarray*}}
\newcommand{\ednn }{\end{eqnarray*}}
\newcommand{\bdmn}{\begin{eqnarray}}
\newcommand{\edmn}{\end{eqnarray}}

\newcounter{biblio}
\newenvironment{references}%
{\begin{list}{[\arabic{biblio}]}{\usecounter{biblio}%
\setlength{\leftmargin}{2.5em}\setlength{\rightmargin}{0pt}%
\setlength{\labelwidth}{2em}\setlength{\itemsep}{0pt}}}{\end{list}}
\newcommand{\References}%
{\vspace{2.8ex plus .3ex minus .3ex}%
\begin{center}{\bf References}\end{center}\begin{references}}

\usepackage{amssymb}
\usepackage{mathrsfs}
\newcommand{\N}{{\mathbb{N}}}
\newcommand{\Z}{{\mathbb{Z}}}
\newcommand{\zd}{\Z^d}

\newcommand{\R}{{\mathbb{R}}}

\newcommand{\tl}{\widetilde}

\newcommand{\limn}{\lim_{n \nearrow \8}}

\newcommand{\vvs}{\vspace{2ex}}
\newcommand{\vs}{\vspace{1ex}}

\newcommand{\lef}{\left}
\newcommand{\rig}{\right}

\newcommand{\8}{\infty}

\newcommand{\dps}{\displaystyle}

\newcommand{\half}{\mbox{$\frac{1}{2}$}}

\renewcommand{\a}{\alpha}
\renewcommand{\b}{\beta}
\newcommand{\gm}{\gamma}

\newcommand{\e}{\varepsilon}
\newcommand{\eps}{\varepsilon}

\newcommand{\h}{\eta}
\newcommand{\tht}{\theta}

\newcommand{\lm}{\lambda}

\newcommand{\m}{\mu}

\newcommand{\s}{\sigma}
\renewcommand{\S}{\Sigma}

\newcommand{\w}{\omega}
\newcommand{\tw}{\tilde \omega}
\newcommand{\om}{\omega}
\newcommand{\W}{\Omega}

\newcommand{\cF }{{\cal F}}
\newcommand{\cG }{{\cal G}}

\newcommand{\cN }{{\cal N}}

\newcommand{\cS }{{\cal S}}

\makeatletter
\def\section{\@startsection{section}{1}{\z@}{-3.5ex plus -1ex minus 
 -.2ex}{2.3ex plus .2ex}{\bf}}
\def\subsection{\@startsection{subsection}{2}{\z@}{-3.25ex plus -1ex minus 
 -.2ex}{1.5ex plus .2ex}{\bf}}
\makeatother

\newcommand{\cvlaw}{\stackrel{\rm{ law}}{\longrightarrow}}

\begin{document}
\parindent=0pt

\bcenter

\large{\bf Weak disorder for low dimensional polymers:\\
The model of stable laws.}

\vvs

\vvs \normalsize

\noindent Francis COMETS\footnote{
 Partially supported by CNRS, UMR 7599 ``Probabilit\'es et Mod\`eles
al\'eatoires'' and
Projet GIP ANR POLINTBIO\\
{\tt http://www.proba.jussieu.fr/pageperso/comets/comets.html}
}
\\

\vs \small
Universit{\'e} Paris 7, \\
Math{\'e}matiques, Case 7012\\
 2 place Jussieu,
75251 Paris, France \\
{\tt email: comets@math.jussieu.fr} \\

\ecenter

\normalsize

\begin{abstract}
In this paper, we consider directed polymers in random environment
with long range jumps in discrete space and time. 
We extend to this case 
some techniques, results and classifications 
known in the usual short range case. However, some properties are 
drastically 
different when the underlying random walk belongs to the domain 
of attraction of an $\a$-stable law. For instance, 
we construct natural examples of directed 
polymers in random environment which experience weak disorder in low 
dimension.

\end{abstract}
\vspace{1cm}
\footnotesize 
\noindent{\bf Short Title.} Stable Directed Polymers
 


\normalsize

\section{Introduction}

Directed polymers in random environment can be viewed as 
random walks in a random potential, which is inhomogeneous 
both in time and space. We restrict here to the discrete case
where the walk has discrete time and space $\Z^d, d \geq 1$.
A number of  motivations for considering
these models are given in the physics litterature, 
in the context of growing random surfaces   \cite{KrSp91}, of 
nonequilibrium steady states and phase transitions \cite{Racz}.
An increasing interest for these models is showing up in the 
mathematical community, and  recent, striking results are the 
characterizations of the ``weak disorder -- strong disorder'' 
and the ``delocalization -- localization'' transitions
given in \cite{CaHu02} and \cite{CSY03}. We give precise 
definitions of these concepts in definition
\ref{def:WD-SD} and above corollary \ref{cor:loc}. Roughly,
weak disorder and delocalization mean that the polymer 
behaves like the random walk, although 
strong disorder and localization mean that it is extremely influenced by 
the medium and it concentrates in just a few corridors where 
the  medium is favorable. It is not known whether these two
phase transitions coincide or not; However, a
partial step is made  for nearest neighbor walks in \cite{CaHu06}.
Heavy tailed  environments are studied in \cite{Va05}, 
they cause a strong form of localization. 
Small dimensions are shown to be special
\cite{CaHu02}, \cite{CSY03}:  for nearest neighbor walks,
strong disorder always in dimension 1 and 2. Moreover, it was
recently proved 
\cite{CoVa05} that the polymer is always localized  in dimension 1.

\medskip

One of the purposes of the present paper
is to clarify the nature of the ``weak disorder -- strong disorder''
transition. We show that 
strong disorder relates to an infinite number of meetings
of two independent random walkers, for a variety of models.
We also study the influence of the jump distribution on 
the ``delocalization -- localization'' transition, and the interplay
between jump tails, space dimension and existence of  delocalized
phase. An interesting contribution here is to  
construct natural, general examples of directed 
polymers in random environment which experience weak disorder in low 
dimension. The jumps have to be long-tailed.
Since the pionneer work of Paul L\'evy, 
the long-time behavior of such walks is known to be classified by stable 
laws, and our results will depend on the stable law which attracts
the random walk.

\medskip

Stable laws and L\'evy flights  model abnormal diffusion
and  mimic rapid turbulent transport.   They also arise naturally from
coarse-graining procedures for short range walks, e.g.  hiting times.
L\'evy flights in a random potential
are considered in \cite{chen-deem} to analyze the $A+A \to \emptyset$
chemical reaction
and explain the phenomenon of superfast reaction, when a small 
amount of potential disorder
added to the turbulent fluid leads to an increase rate of the reaction.
Dynamics of particle randomly moving along a disordered hetero-polymer
subject to rapid conformal changes, lead to superdiffusive motion. 
A model is introduced in
\cite{brockmann-geisel}, corresponding to a L\'evy flight in a 
random potential in chemical coordinates. Both these models have 
time-independent potential, but the time-dependent case simply corresponds
to crossings in the presence of strong external fields.
An instructive review of the ocurrence of L\'evy processes in 
sciences as fluid mechanics, solid state physics, polymer chemistry
and mathematical finance, is given in \cite{woy}.

\medskip

We will assume that  the  random walk belongs to the 
domain of attraction of an $\a$-stable law for some $\a \in (0,2]$. 
This implies that random walk at large times $n$ roughly 
scales like $n^{1/\a}$. In  the case $\a <2$
it also means that the tails $P(|\w_1|>r)$ of individual jumps
are of order $r^{-1/\a}$ for large $r$. The case $\a=2$ includes the usual
one where the walk is nearest neighbor.
The medium is assumed to have finite  exponential moments. 
We prove that weak disorder holds for $d=1, \a<1$ and $d=2, \a <2$, at least 
for high temperature. 
This is rather surprising in view of the results mentioned above in
the simple random walk case.
In dimension $d\geq 3$, our results here are not qualitatively different from 
those obtained for the simple random walk. For completeness we will 
state the results in all dimensions, but we emphasize low dimensions.
All through, we assume that the environment has finite exponential moments.
\medskip

The paper is organized as follows. In the next section, we introduce the model
and recall some necessary facts on stable laws and their attraction domain.
Then, the free energy is defined, together with the regimes of
weak disorder and of strong disorder. We give sufficient conditions 
for weak disorder in section \ref{sec:3} together 
with some properties of the polymer 
there, and  sufficient conditions 
for strong disorder in section \ref{sec:4}. 
The last section is dedicated to the phase diagrams and localization 
properties.
As already mentioned, we extend  some constructions and 
results from nearest neighbor random walks to long range ones,
we will not repeat proofs unless necessary but indicate precise references 
instead.

\section{Long jumps polymers} \label{sec:discrete}

We first need to state a few elementary facts on 

\subsection{Stable laws} 

These are all possible distributional limit of sums of i.i.d.
random vectors up to renormalization. 
By definition, a stable law on $\R^d$ is such that, for all $n \geq 1$,
if $X_1,\ldots X_n$ are i.i.d. with this law, there exist $a_n>0$ and 
$b_n \in \R^d$ such that 
$$
 \frac{X_1+\ldots+ X_n - b_n }{a_n}\quad {\rm still\ has\ this \ law}
$$ 
To avoid triviality we assume that the law is not a Dirac mass.
It can be shown that there is a 
unique $\a \in (0,2]$ such that for all $n$, the above $a_n$ is 
$a_n=n^{1/\a}$. This exponent $\a$ is called the {\it index} 
(or characteristic exponent) of the stable law,
and we also say that $P$ is $\a$-stable. 

Except a few special cases, stable laws are complicated, 
they cannot be written in terms of simple 
functions, but their Fourier transforms  are simple. An $\a$-stable 
random vector $S_\a$ has characteristic functions
$$E(e^{iz\cdot S_\a})=e^{\psi(z)}\;,\quad 
z \in \R^d$$ 
where the form of the exponent $\psi$ depends on the  index $\a \in (0,2)$:
\begin{equation} \label{eq:caracst2}
{\rm for} \ \a=2\;, \quad \psi(z)=
 i \tau \cdot z - \half z \cdot Az
\end{equation}
with $\tau \in \R^d$ and $A$ a $d\times d$ symmetric positive definite
matrix;
\begin{eqnarray}  \label{eq:caracsta}
{\rm for} \ \a\neq 1,2\;, \quad \psi(z)&=&\psi(z)_{\a,\tau,\sigma}
\\  \nonumber
&=& i \tau \cdot z - 
\int_{\cS^{d-1}} |z\cdot \xi|^\a  \left(1- i \tan \frac{\pi \a}{2}\;
{\rm sgn}( z\cdot \xi ) \right) \sigma( d \xi)  
\end{eqnarray}
with $\tau \in \R^d$ and $\sigma$ a finite nonzero measure on the unit 
sphere $\cS^{d-1}$ (the sign function is defined by sgn$(u)=1$ for
$u>0$,  sgn$(u)=-1$ for
$u<0$ and sgn(0)=0);
\begin{eqnarray}
\label{eq:caracst1}
{\rm for} \ \a= 1\;, \quad \psi(z)&=&\psi(z)_{\a,\tau,\sigma}\\
\nonumber &=& i \tau \cdot z - 
\int_{\cS^{d-1}} \left( |z\cdot \xi| + \frac{2i}{\pi} z\cdot \xi 
\log |z\cdot \xi| \right) \sigma( d \xi)  
\end{eqnarray}
with $\tau \in \R^d$ and $\sigma$ a finite nonzero measure on the unit 
sphere $\cS^{d-1}$. The vector $\tau$ 
(sometimes called the translate)
and the measure $\sigma$ (sometimes called the spherical part of the
L\'evy measure)
are uniquely defined. They are location and asymmetry 
parameters. The law is invariant under rotations centered at some $x \in \R^d$
if and only if $x=\tau$ and $\sigma$ a uniform measure on the sphere; 
The law is invariant under the central symmetry with center $x \in \R^d$
if and only if $x=\tau$ and $\sigma$ is invariant under $\xi \mapsto -\xi$.

Here are some special cases where the density is simple.  
In the case $\a=2$, the law is the $d$-dimensional Gaussian 
with mean $\tau$
and covariance matrix $A$, with density $$ x \mapsto
(2\pi)^{-d/2} ({\rm
  det}A)^{-1/2} \exp\left\{ -\frac{1}{2} (x-m)^* A^{-1}(x-m)\right\}\;.$$


For $c>0$, $\tau \in \R^d$ and $\Gamma$ the Euler function,
the $d$-dimensional Cauchy law with density 
$$   x \mapsto \Gamma((d+1)/2) 
\frac{c } 
{\pi^{(d+1)/2}  (|x-\tau|^2+c^2)^{(d+1)/2}}
$$
is stable with $\a=1$, with characteristic exponent
$\psi(z)_{1,\tau,\sigma}$ with $\sigma$ the uniform measure of mass 
$c$.

A complete overview on stable laws and domains of attraction 
is given in the book \cite{BGT}, and a shorter presentation in 
\cite{brei}. For  stable processes, 
we refer to  the books \cite{BGT} and \cite{sato}.

\subsection{The model}

$\bullet${\it The random walk:} 
$(\{ \w_n\}_{n \geq 0}, P)$ is a random walk on $\Z^d$ starting from 0, 
ie, the variables $\w_{k+1}-\w_k (k=1,2, \dots)$ are i.i.d. under $P$ with
$\w_0=0$, and
we denote by $q$ their common law $q(x):=P(\w_1=x)$.
We assume that $q$ belongs to the  domain of 
attraction of a stable law (on $\R^d$) with some index $\a \in (0,2)$. More 
precisely, we assume that there exist $\a \in (0,2], \tau \in \R^d,
\sigma$ a finite nonzero measure on $\cS^{d-1}$, and deterministic  sequences
$a_n>0, b_n \in \R^d$,  such that
\begin{equation}\label{hyp:domst}
P\left( \exp\{i z \cdot  \frac{\w_n -b_n}{a_n}\} \right) 
\longrightarrow \; \exp \psi_{\a,0, \sigma}(z)
\end{equation}
for all $z$. To simplify our discussion, we will also assume that 
the limit is truly $d$-dimensional, ie, that it satisfies (\ref{hyp:2d})
below.

We now give a short account on  our assumption  (\ref{hyp:domst}),
and recall some facts on  the  domain of 
attraction of stable laws, starting with the simpler case of dimension $d=1$. 
The reader may also decide to skip these details in a first reading, and jump 
directly to the important example \ref{ex:attr1}.

\medskip

{\bf Dimension 1:} In one dimension, 
this assumption can be described in terms of  
the tails of $q$. 
We follow the presentation of  section 8.3 in \cite{BGT}.
The  cases $\a \in (0,2)$ and $\a=2$ being different, we start with

\begin{enumerate}
\item Case $\a \in (0,2)$. We let
$R_0$ be the space of slowly varying functions in the sense of Karamata, 
i.e. of functions
$\ell: [0, \8)\mapsto [0, \8)$ such that 
$$
\ell(sr)/\ell(r) \to 1 \quad (r \to \8) \quad \forall  s>0.
$$
Examples of such functions 
are constants, $\ln r$ or $\exp [\ln r / \ln \ln r]$.
Assumption (\ref{hyp:domst})
is equivalent to 
\begin{equation}\label{eq:domst1}
P(|\w_1|\geq r) = r^{-\a} \ell(r)
\end{equation}
for some $\ell \in R_0$, and 
\begin{equation}\label{eq:domst2}
\frac{P(\w_1\leq -r)}{P(|\w_1|\geq r)} \to q_*\;,\quad
\frac{P(\w_1\geq r)}{P(|\w_1|\geq r)} \to p_*\quad(r \to \8)
\end{equation}
where we note that $p_*+q_*=1$. Then, the sequence $a_n$ is of the form
$$
a_n=n^{1/\a} \ell'(n)
$$
with a  slowly varying function $\ell' \in R_0$ which can be taken 
such that 
$$P^2(e^{i z (\w_1-\tw_1)/a_n})^n \to e^{- 2|z|^\a} \quad (n \to \8)\quad \forall
z \in \R$$
In this case, the limit in  (\ref{hyp:domst}) has 
$$\sigma(+1)=p_*\;,\quad
\sigma(-1)=q_*$$
(note that $S^0=\{-1, +1\}$).  Moreover, we can take 
$$b_n=0 \quad {\rm for} \; \a <1\;,$$ 
$$b_n=n P(\w_1) \quad {\rm for} \;\a \in (1,2)\;,$$ 
in which case $\w_1$
is integrable, and 
$$b_n=n P[\psi(\w_1/a_n)]\;,\quad 
\psi(t)=t(1+t^2)^{-1} \quad {\rm for} \; \a =1\;.$$ 
The reader is refered to  \cite{BGT}, pp. 343-347,
for further details.

\item $\a=2$: the assumption (\ref{hyp:domst})
is then equivalent to the function $r \mapsto P(|\w_1|^2: |\w_1|\leq r)$
being  slowly varying. A sufficient condition is 
existence of a second moment of $\w_1$, and (\ref{hyp:domst}) is simply the
standard central limit theorem. 
\end{enumerate}

{\bf Dimension $d \geq 2$:} Let $\varphi(z)=P(e^{iz\cdot \w_1})$.
The characterization of assumption (\ref{hyp:domst}) is known in terms
of the characteristic function $\varphi$ of the law $q$ and of slowly 
varying functions
(see  theorem 2.6.5 in
\cite{ibralin} and \cite{aarden},
corollary 1--2 in section 2).
Assumption (\ref{hyp:domst}) with a truly $d$-dimensional limit 
in the sense of (\ref{hyp:2d}) is equivalent, for $\a \in
(0,2)\setminus \{1\}$, to
$$
\log \varphi(z) = \left\{
\begin{array}{cc}
\phi_{\a,0,\sigma}(z) \ell(1/|z|) + i z \cdot \tau + o(|z|^\a
\ell(1/|z|))
&   {\rm if} \; \a >1\;,\\
\phi_{\a,0,\sigma}(z) \ell(1/|z|)  + o(|z|^\a
\ell(1/|z|))
&   {\rm if} \; \a <1\;.
\end{array}
\right.
$$
This is for $\a \neq 1$, the case $\a=1$ being more complicated.
 We simply mention that, for a {\it symmetric} law $q$ and  $\a=1$,
 assumption (\ref{hyp:domst}) with (\ref{hyp:2d}),  is equivalent to 
$$
\ln \varphi(z) = -
\int_{\cS^{d-1}}  |z\cdot \xi| \sigma( d \xi) \times
 \ell(1/|z|)  + o(|z|
\ell(1/|z|))\;,
$$
and we refer to 
corollary 2 in section 2 of \cite{aarden}, for the somewhat cumbersome
general case.

We give a generic example
where the assumption holds.

\begin{example} \label{ex:attr1} 
Let $d=1$ or $2$,  and $q$ be symmetric with
\begin{equation} \label{eq:jumps1.1}    
q(x) = b(|x|^{-d - \a}+\eps(x))
 \quad {\rm as}\; |x| \to \8\;, \; x \in Z^d
\end{equation}
with $\a \in (0,2)$.
Then, the assumption  (\ref{hyp:domst}) hold true with 
$\tau=0$ and $\sigma$ uniform,
using
\cite{aarden},\cite{ibralin}), and estimates  of the characteristic
function (see e.g. (2.13) and (3.11)
in \cite{giwe}).
\end{example}

\medskip

$\bullet$ {\it The random environment:}  
$\h =\{\h (n,x) : n \in \N,\; x \in \zd \}$ is a 
sequence of r.v.'s which are real valued,
non-constant, and i.i.d.(independent identically distributed) r.v.'s 
 defined on a probability space $(\W_\h, \cG, Q)$ such that 
$$ 
Q[\exp (\b \h (n,x))] <\8 \; \; \; \mbox{for all $\b \in \R$.}
$$
We then let $\lm(\b)=\ln Q[\exp (\b \h (n,x))]$.
\medskip

$\bullet$ {\it The polymer measure:}
For any $n>0$, define the probability measure $\m_n$ on the path space 
$(\W_\om, \cF)$ by  
\begin{equation} \label{mnen}
\m_n (d\w )=\frac{1}{Z_n}
\exp\{ \beta H_n(\om)\} \;P(d\w),
\end{equation}
where $\b >0$ is a parameter (the inverse temperature), where
\begin{equation} \label{Ham}
H_n(\om)=H_n(\h,\om)=\sum_{1 \leq j \leq n} \h (j,\w_j)
\end{equation}
 and 
\begin{equation} \label{Zn}
Z_n=Z_n(\beta, \h)
=P\lef[\exp \lef(\b\sum_{1 \leq j \leq n} \h (j,\w_j)\rig) \rig]
\end{equation}
is the the partition function.

\section{Free energy, and the natural martingale} \label{sec:free}

The partition function is random, but it is self-averaging as $n$ increases.

\begin{proposition} \label{prop:freeenergy} Let $\a \in (0,2]$ arbitrary.
  As $n \to \8$, the quenched free energy converges 
to a deterministic constant:
\begin{equation} \label{eq:freeenergy}
 \frac{1}{n} \ln Z_n \longrightarrow  p(\beta) := \lim_{m \to \8}
\frac{1}{m} Q \ln Z_m 
\end{equation}
$ Q-$a.s. and in $L^q$ ($1 \leq q < \8$). Moreover, we have the annealed bound 
\begin{equation} \label{eq:annealedbound}
p(\b) \leq \lm(\b)
\end{equation}
\end{proposition}
$\Box$ The proof in the case of a simple random walk $P$
(proof of prop. 2.5 in \cite{CSY03}, pp. 720--722)
covers the general case of $\a \in (0,2]$
without change. The last inequality comes from Jensen 
inequality, which writes
$$\frac{1}{m} Q \ln Z_m \leq \frac{1}{m}  \ln Q Z_m = \lm (\b)$$ 
\qed
The sequence $(W_n, 
n \geq 1)$ defined by
\begin{equation} \label{Wn}
W_n=Z_n \exp (- n \lm (\b )) 
\end{equation}
is   a  positive, mean 1, martingale with respect to the environmental 
filtration
$(\cG_n)=\sigma(\h(t,x), t \leq n, x \in \Z)$. This was noticed first by 
Bolthausen \cite{Bol89}. 
By the martingale convergence theorem, the limit $$W_\8=\limn W_n$$ exists 
$Q$-a.s. It is clear that 
the event $\{ W_\8=0\}$ is measurable with respect to the tail $\s$-field
$
\bigcap_{n \geq 1}\s [ \h (j,x) \; ; \; j \geq n, \; x \in \zd ]\;.
$
By Kolmogorov's zero-one law 
 every event in  the tail 
$\s$-field has probability 0 or 1. Hence,
 there are only two possibilities for the 
positivity of the limit
\begin{equation} \label{Z8>0}
Q\{W_\8>0 \}=1\;,
\end{equation}
or 
\begin{equation} \label{Z8=0}
Q\{W_\8=0 \}=1\;.
\end{equation}
\begin{definition} \label{def:WD-SD}
The above situations 
(\ref{Z8>0}) and (\ref{Z8=0}) will be called  the {\bf weak disorder}
phase and the {\bf strong disorder} phase, respectively. In the
first case, $p=\lm$.
\end{definition}

In corollary \ref{cor:loc0} at the end of the paper
 we explain why it is important to decide if the inequality
in (\ref{eq:annealedbound}) is an equality or not.

\section{Existence of weak disorder, properties} \label{sec:3}

We need now to consider on the product space $(\W^2, \cF^{\otimes 2})$,
the probability measure $P^{\otimes 2}=P^{\otimes 2}(d\w,
d\tl{\w})$, 
that we will view as the distribution of the couple $(\w, \tl{\w})$ with 
 $\tl{\w}=(\tl{\w}_k)_{k \geq 0}$  an independent copy 
of $\w=(\w_k)_{k \geq 0}$. 

When $P$ satisfies (\ref{eq:jumps1.1}), we see that the random walk
$\w-\tw$ is attracted by the symmetric $\a$-stable law.
Precisely, with $a_n$ from (\ref{eq:jumps1.1}), we have
\begin{equation} \label{eq:111}
P^{\otimes 2}\left( \exp\{i z \cdot  \frac{\w_n-\tw_n}{a_n}\} \right) 
\longrightarrow \; \exp \left( \psi_{\a,0, \sigma}(z)+\psi_{\a,0,
    \sigma}(z)
\right)=
\exp \psi_{\a,0, \sigma'}(z)
\end{equation}
with $\sigma'(B)=\sigma(B)+\sigma(-B)$ 
for all Borel subset $B$ of 
$\cS^{d-1}$.

For later purposes, it is essential to observe that the difference 
$\w-\tw$ is a transient
random walk 
-- i.e., $ N_\8:=\sum_{n=1}^\8 {\bf 1}_{\w_n=\tw_n}<\8$ a.s. -- 
in the three following cases:
\begin{description}
\item  (i) $d=1$ and $\a \in (0,1)$,
\item (ii) $d=2$ and $\a \neq 2$,
\item (iii) $d \geq 3$ and $\a \in (0,2]$,
\end{description}
provided the limit is truly
$d$-dimensional. This extra assumption for $\a \in (0,2)$ means 
that the linear space
spanned by the support of the measure $\sigma$ is $\R^d$, and 
for $\a=2$ that the covariance matrix $A$ is non-degenerate,
\begin{equation} \label{hyp:2d}
{\rm Vect} (\, {\rm supp}\,
 \sigma ) = \R^d \quad {\rm or}  \quad {\rm rank}(A)=d \;,
\end{equation}
according to the case $\a<2$ or $\a=2$.
Indeed, with $\phi(z)=P^{\otimes 2}[ \exp\{i z \cdot (\w_1-\tw_1)\} ]$, 
(\ref{eq:111}) amounts to 
\begin{equation} \label{eq:112a}
\phi(z) =  
\exp\Big\{\psi_{\a,0, \sigma'}(z) \ell\Big(\frac{z}{|z|},\frac{1}{|z|}
\Big)\Big\}
\end{equation}
with $\ell(\xi,\cdot)$ a slowly varying function depending continuously
on $\xi \in {\mathcal S}^{d-1}$. 
Since it holds, under   (\ref{hyp:2d}),
$$\psi_{\a,0, \sigma'}(z) \leq C 
|z|^\a\;,$$ we have
\begin{equation} \label{eq:112}
\int_{[-\pi, \pi[^d} \frac{dz}{1-\phi(z)}
\quad 
\left\{
\begin{array}{lll}
=\8 &{\rm if}& \a > d\\  
<\8 &{\rm if}& \a < d
\end{array}
\right.  
\end{equation}

Applying the Chung-Fuchs criterion (P1 in section 8 of \cite{Sp}),
we see that the walk $\w-\tw$ is transient in the second case, and then
$$
\pi(p):=P^{\otimes 2}
(\exists n \geq 1 \; \w_n-\tw_n = 0) <1
$$ 

\begin{remark}
The walk $\w-\tw$ is recurrent when $ \a > d$. The  border case $\a=d$
is more  subtle:
Transience may hold or may not hold 
depends on the 
 slowly varying term in (\ref{eq:112a}). In the positive, 
 the validity of the
the next two theorems will extend to critical cases $\a=d=1$, $\a=d=2$.
\end{remark}

\begin{theorem} \label{th:WDa} {\bf Weak disorder region in dimension 1, 2.}
In addition to (\ref{hyp:2d}), 
assume  either (i) $d=1$ and $\a \in (0,1)$, or (ii) $d=2$, $\a \neq
2$, or (iii) $d \geq 3$ and $\a \in (0,2]$. 
Then, for all $\b$ such that 
\begin{equation} \label{L2}
\lm(2\b)-2\lm(\b)<  \ln 1/\pi(p)\;,
\end{equation}
 we have $W_\8 >0$\ $Q$-a.s.
\end{theorem}

The result may come as a  surprise, since for $P$ the simple random walk, 
it was proved that 
 $W_\8 =0$ $Q$-a.s. \cite{CaHu02} \cite{CSY03}, and even that $p<\lm$
\cite{CoVa05}. The method used in the last reference is based on comparisons 
with polymers models on trees. It is impossible to extend it to 
long range jumps, although related
ideas can be --and will be-- in the sequel, see (\ref{eq:KP}). 

Following the techniques of \cite{Birk} using a conditional second
moment, one could extend the validity of the result to a domain in
$\b$
larger than (\ref{L2}).
 
\medskip

$\Box$ Following  \cite{Bol89}
we compute the $L^2$-norm 
of the martingale $W_n$. To do so, 
 we represent $W_n^2$ in terms of an independent couple $(\w, \tl{\w})$
introduced above.
\begin{eqnarray*}
Q[W_n^2]&= &
Q\left[ P^{\otimes 2} \prod_{t=1}^n e^{\b [\h(t,\w_t)+
\h(t,\tw_t)]-2 \lm(\b)}\right]\\
&=&  P^{\otimes 2}\left[ \prod_{t=1}^n \left(e^{\lm(2 \b)-2 \lm(\b)}
 {\bf 1}_{\w_t=\tw_t} +  {\bf 1}_{\w_t \neq \tw_t}\right)\right] \\
&=&  P^{\otimes 2} \left[ e^{ \gamma_1 N_n} \right]\;,
\end{eqnarray*}
with $\gamma_1=\lm(2\b)-\lm(\b)$, and 
$N_n$ the number of intersections of the paths $\w, \tw$ up to
time $n$,
\begin{equation} \label{Nn}
N_n=N_n(\om, \tw) = \sum_{t=1}^n {\bf 1}_{\w_t=\tw_t}
\end{equation}
As $n \to \8$, $N_n \nearrow N_\8$, and by monotone convergence
$Q[W_n^2] \nearrow P^{\otimes 2} \left[ e^{C N_\8}\right] $.

In the cases under consideration, the random variable $ N_\8$ 
is geometrically distributed 
$$
P^{\otimes 2}(N_\8=n)=\pi(\sigma)^n[1-\pi(\sigma)]\;,\quad n\geq 0\;,
$$
with $\pi(p)$ the probability of  return defined above the theorem.
Hence it has finite exponential moments
$$
P^{\otimes 2}[ \exp\gamma_1  N_\8 ] <\8 \; \iff \; \gamma_1 < \ln 1/\pi(p)
$$ 
Therefore, when $\lm(2\b)-2\lm(\b)<  \ln 1/\pi(p)$, the martingale $W_n$ is
bounded in $L^2$, and by the classical $L^2$-convergence theorem, 
it converges
in $L^2$ to a limit, which is necessarily equal to $W_\8$. So
$Q  W_\8 = \lim_n Q  W_n = 1$, which excludes the possibility that 
the limit vanishes.
 \qed
\medskip

Inside the subset of the weak disorder region determined by the condition 
(\ref{L2}), the fluctuations of the path 
remain similar to those of $P$.

\begin{theorem} \label{th:WDaf} 
Assume  either (i) $d=1$ and $\a \in (0,1)$, or (ii) $d=2$, $\a \neq
2$, or (iii) $d=3$ and  $\a \in (0,2]$, 
and assume (\ref{hyp:2d}). When (\ref{L2}) holds,
we have for all 
bounded continuous function $g: \R \to \R$,
$$
\mu_n \left[ g \left( \frac{\w_n-b_n}{a_n}\right) \right]
\to \nu( g )
$$
in $Q$-probability as $n \to \8$, where $\nu$ is the $\alpha$-stable law
with characteristic function $\psi_{\a,0,\sigma}$. 
\end{theorem}

$\Box$ We let $\nu_n(\cdot) = P[ (\w_n-b_n)/a_n \in \cdot]$. 
\begin{eqnarray} \nonumber
 && Q \left( \left\vert \mu_n \left[ g\left( \frac{\w_n-b_n}{a_n}\right)
      \right] 
- \nu_n(g) \right\vert^2 W_n^2\right)
\\ \nonumber
&=& 
P^2  Q  \left( e^{\b H_n(\w)+\b H_n(\tw)-2n\lm}
\left[  g\left( \frac{\w_n-b_n}{a_n}
      \right)
- \nu_n(g)\right]
\left[  g\left( \frac{\tw_n-b_n}{a_n}\right)
- \nu_n(g)\right]
 \right)
\\ \label{eq:100}
&=&
P^2   \left( e^{\gamma_1 N_n}
\left[  g\left( \frac{\w_n-b_n}{a_n}
      \right)
- \nu_n(g)\right]
\left[  g\left( \frac{\tw_n-b_n}{a_n}\right)
- \nu_n(g)\right]
\right)
\end{eqnarray}
We know that, under $P^2$,  the r.v. $N_n$ converges to $N_\8$ a.s., and that 
$(\w_n-b_n)/a_n$ -- and similarly 
$(\tw_n-b_n)/a_n$ -- converges to $\nu$
in law. 
Now, we claim that,  under $P^2$, the triple 
\begin{equation}\label{eq:101}
(N_n,(\w_n-b_n)/a_n,
(\tw_n-b_n)/a_n) \cvlaw
(N,S,\tilde S)
\end{equation}
with $(N,S,\tilde S)$
an independent triple 
where $N$ has the same law   
as $N_\8$, $S$ and $\tilde S$ have the law $\nu$.  The proof of this fact
makes use of the observation that 
$$\sup_{n \geq m} P^2(N_n \neq N_m) \to 0 \; {\rm as } \; n \to \8$$
since $N_n \nearrow N_\8 < \8$ a.s.
Fix $m \geq 1$ and  $f,g,\tilde g$ continuous and bounded. For all $n \geq m$,
we write
\begin{eqnarray*}
\!\!\!\!\!\!\!\!\!\!\!\!\!\!\!\!\!\!\!\!\!
&& 
\!\!\!\!\!\!\!\!\!\!\!\!\!\!\!\!\!\!\!\!\!
P^2   \left[ f(N_n) g\big( \frac{\w_n-b_n}{a_n}
\big) \tilde
  g\big( \frac{\w_n-b_n}{a_n}
\big)
\right]
\qquad \qquad \qquad \qquad \qquad \qquad \qquad 
 \\
&=&
P^2   \left[ f(N_n) g\big( \frac{\w_n-b_n}{a_n}
\big) \tilde
  g\big( \frac{\w_n-b_n}{a_n}
\big)  {\bf 1}_{N_n = N_m}
\right]
 + \eps(n,m)\\
&=&
P^2   \left[ f(N_m) g\big( \frac{\w_n-b_n}{a_n}
\big) \tilde
  g\big( \frac{\w_n-b_n}{a_n}
\big)  {\bf 1}_{N_n = N_m}
\right]
 + \eps(n,m)
\\
&=&
P^2   \left[ f(N_m) g\big( \frac{\w_n-\w_m-b_n}{a_n}
\big) \tilde
  g\big( \frac{\w_n-\w_m-b_n}{a_n}
\big)  {\bf 1}_{N_n = N_m}
\right]
 + \eps'(n,m)\\
&=&
P^2   \left[ f(N_m) g\big( \frac{\w_n-\w_m-b_n}{a_n}
\big) \tilde
  g\big( \frac{\w_n-\w_m-b_n}{a_n}
\big) 
\right]
+ \eps"(n,m)
\\
&=&
P^2   [ f(N_m)] \times P\left[ g\big( \frac{\w_n-\w_m-b_n}{a_n}
\big)\right] \times 
 P\left[  \tilde g\big( \frac{\w_n-\w_m-b_n}{a_n}
\big) 
\right]
+ \eps"(n,m)\;,
\end{eqnarray*}
which equalities define the terms $\eps(n,m), \eps'(n,m), \eps''(n,m)$
on their first occurence.
Here, 
$$|\eps(n,m)| \leq \|f\|_\8  \|g\|_\8  \|\tilde g \|_\8  
P(N_n \neq N_m)$$ tends to 0 as $m \to \8$ uniformly in $n \geq m$,
$\eps'(n,m)-\eps(n,m) \to 0$ as $n \to \8$ for all fixed $m$, 
and $\sup_{n \geq m}\eps''(n,m) \to 0$ as $m \to \8$. The last equality
comes from independence in the increments of the random walks, and of the
two random walks $\w$ and $\tw$. Hence, letting $n \to \8$ and then  
$m \to \8$, we get
$$
P^2   \left[ f(N_n) g\big( \frac{\w_n-b_n}{a_n}
\big) \tilde
  g\big( \frac{\w_n-b_n}{a_n}
\big)
\right]
\to 
P^2   [ f(N_\8)] \times \nu[ g] \times  \nu[\tilde g]
$$
which proves (\ref{eq:101}). Coming back to (\ref{eq:100}),
and since $P^2( e^{\gamma N_n})<\8$ for some  small enough
$\gamma>\gamma_1$, (\ref{eq:101}) implies that
$$ Q \left( \left\vert \mu_n \left[ g\left( \frac{\w_n-b_n}{a_n}\right)
      \right] 
- \nu_n(g) \right\vert^2 W_n^2\right)
\to 
P^2   ( e^{\gamma_1 N_\8}) \left[
\nu (g)
- \nu(g)\right]^2
=0
$$
Since $W_n^{-2}$ converges to a finite limit, it is bounded in probability,
this yields the desired convergence in probability.
\qed

\begin{remark} \label{rem1}
(i) In the case when $P$ is the nearest neighbor simple random walk,
the condition (\ref{L2}) implies a quenched central limit 
theorem, ie, 
that  central limit 
theorem holds for a.e. realization of the environment
 \cite{Bol89}. Our result here is weaker. Due to the lack of moments
for the long jumps here, the natural martingales which can be used
in the standard case are not defined in the present setup.\\
(ii) In the case when $P$ is the nearest neighbor simple random walk,
it was shown in \cite{CY05} that central limit theorem holds (in a weak form
at least) as soon as $W_\8>0$. Then it is questionable whether in the 
model of the 
present paper, weak disorder implies  convergence of the renormalized
position of the polymer to an $\a$-stable law. We leave the question open.\\
(iii) When $P$ is the simple random walk, many other results are known 
under condition (\ref{L2}), for instance:
\begin{enumerate}
\item 
Local limit theorem for 
the polymer measure \cite{Sin95},  \cite{Va04};
\item  How does the polymer 
depends on the environment ? (This question is answered in \cite{BMP}
by computing the random corrections to gaussian for 
cumulants of the polymer position.)
\end{enumerate}
We leave open the question of which is the  counterpart of these results
for long range random walks we consider here ($\a \in (0,2)$).
\end{remark}


We end this section with a model where weak disorder holds 
at all temperature and all dimension. Viewed as a growing random
surface, it does not have a roughening transition, and consequently
in this respect,
it does not belong to the Kardar-Parisi-Zhang (KPZ) class
\cite{KrSp91}.

\begin{example}\label{CLTB}
{\it Bernoulli environment.}
The case when $\h(t,x)=1$ or 0
with probability $q$ and $1-q$ respectively, is remarkable since 
weak disorder may hold at all temperature. 
Here we find $\lm= \ln [ p e^{\b}+(1-q)]$, and we see from direct
computations that 
$${\dps 
\lim_{\b \nearrow \8}\gm_1 (\b)=
-\ln(q).}$$
Hence, (\ref{L2}) holds for all $\b \geq 0$ if 
 $q > \pi(p)$.
Theorem \ref{th:WDa} shows that, in this case, weak disorder
holds for all $\b \geq 0$, and Theorem \ref{th:WDaf} shows that
the polymer position at time $n$ still fluctuates at order $n^{1/\a}$.
\end{example}

\section{Existence of strong disorder} \label{sec:4}

The next result gives a sufficient condition for $p<\lm$, which implies 
 strong disorder.

\begin{proposition} \label{th:WD} Let  $\a \in (0,2]$ and $d$ arbitrary. 
If
\begin{equation} \label{eq:KP}
\b \lm'(\b)-\lm(\b) > -\sum_{k \in \Z}  q(k) \ln q(k)
\end{equation}
then $p<\lm$.
\end{proposition}

We note that the important quantity is here the entropy 
$- \sum_x  q(x) \ln q(x)$ of the walk, which does not directly relates
to the recurrence/transience behavior of the walk.
The entropy is always finite under our assumptions on $q$. 
\begin{example}\label{CLTG}
{\it Gaussian environment.} If $\eta$ is standard gaussian $\cN (0,1)$,
then $\gm_1 (\b)=\b^2$ and hence (\ref{L2}) holds if 
$\b <\sqrt{\ln (1/\pi_d)}$, though (\ref{eq:KP}) holds for 
all $\b > - \sum_x  q(x) \ln q(x)$. In this case, a phase transition takes 
place between weak and strong disorder.
\end{example}

$\Box$ 
Note that
\bdnl{Znm}
Z_{n}= \sum_{x \in \Z} q(x) e^{\b \h(1,x)} Z_{1,n}^x
\;, 
\edn 
where $ Z_{1,n}(x)$ has the same law as $Z_{n-1}$.
Let $\tht \in (0,1)$. By the subadditive estimate 
$$(u+v)^\tht \leq u^\tht +v^\tht\;,\quad 
u, v >0\;,$$ 
we get
$$
Z_n^\tht \leq 
\sum_{x \in \Z}  q(x)^\tht e^{\b\tht \h(1,x)} (Z_{1,n}^x)^\tht
$$
Since $Z_{1,n}^x$ has the same law as $Z_{n-1}$, we obtain a bound on 
$u_n=Q Z_n^\tht$:
\begin{eqnarray*} 
u_n &\leq& \left[ \sum_{x \in \Z}  q(x)^\tht \exp\{\lm( \b \tht)\}\right]
u_{n-1}\\
&\leq& \left[ \sum_{x \in \Z}  q(x)^\tht \exp\{\lm( \b \tht)\}\right]^n
\end{eqnarray*}
by induction. Now, observe that
$$
Q \frac{1}{n} \ln Z_n =
Q \frac{1}{n\tht} \ln Z_n^\tht \leq
  \frac{1}{n\tht} \ln Q Z_n^\tht
$$
which, combined with the previous bound on $u_n=Q Z_n^\tht$, yields
$$
p \leq \inf_{\tht \in (0,1)} \left\{ \frac{1}{\tht} v(\tht) \right\}\;,\quad
v(\tht)=
 \left[ \lm( \b \tht) + \ln
\sum_{x \in \Z}  q(x)^\tht \right]\;.
$$
Since the function $v$ is convex and  positive at 0,
there are only possibilities for the infimum. If the derivative
of $v(\tht)/\tht$ is positive at $\tht=1$, the infimum is acheived at some 
$\tht \in (0,1)$, and is strictly less than $\lm(\b)$ (which is the value at 
1). On the contrary, if  the derivative
of $v(\tht)/\tht$ is less or equal to 0, the infimum is for $\tht \to 1$
and the value is $\lm(\b)$.
Finally, we compute
$$\frac{d}{d\tht} \frac{v(\tht)}{\tht}_{|\tht=1} = \b \lm'(\b)-\lm(\b) +
\sum_{k \in \Z}  q(k) \ln q(k) >0$$
which proves the claim.
 \qed

\section{Phase diagram, transitions and localization}

The sum in (\ref{eq:KP}) is always finite with our choice of $q$. 
For  unbounded $\h$'s, one easily checks that $\lim \b \lm'(\b) -  \lm(\b)
=+ \8$ as $\b \to +\8$, so that  
this condition (\ref{eq:KP}) will be checked for large $\b$. On the other 
hand, the set of $\b$'s such that $p=\lm$
is an interval (possibly with length 0). Indeed it is readily checked that 
the argument in Th. 3.2.(b) in \cite{CY05} for the case $\a=2$,
extends to all values of $\a \in (0,2]$.
To summarize,
\begin{theorem} \label{th:phdi}  {\bf Phase diagram.}
The exists a $\b_c \in [0,\8)$ such that 
$p=\lm$ for $\b \in [0, \b_c]$ and $p< \lm$ for $\b > \b_c$.
Moreover, under the assumption   (\ref{hyp:2d}), we have
$$\b_c >0$$ 
in the following cases:
 (i) $d=1$ and $\a \in (0,1)$, or (ii) $d=2$, $\a \neq
2$, or (iii) $d=3$ and  $\a \in (0,2]$. 
\end{theorem}

When $\a=2$ it is well known  \cite{CaHu02} \cite{CSY03}
that the discrepancy between $p$ and $\lm$ 
relates to localization property of the polymer. In fact the computations
in the special case $\a=2$ still work for all $\a$ (e.g.
\cite{CSY03},
th. 2.1 and proof pp. 711--715). 
Then, we have
\begin{theorem}\label{th:equiv2}
Let $\b \neq 0$, $\a$ and $d$ arbitrary. Define
$$
I_n=\m_{n-1}^{\otimes 2}(\w_n =\tl{\w}_n)$$
Then,
\begin{equation} \label{equiv1}
\{ W_\8 =0 \} = \Big\{ \sum_{n \geq 1}I_n =\8 \Big\}, \; \; \; 
\mbox{$Q$-a.s.}
\end{equation}
Moreover, if $Q\{ W_\8 =0 \}=1$,
there exist $c_1,c_2 \in (0,\8)$ such that $Q$-a.s.,
\begin{equation} \label{equiv2}
 c_1\sum_{1 \leq k \leq n}I_k
\leq -\ln W_n \leq 
 c_2\sum_{1 \leq k \leq n}I_k
\qquad \mbox{for large enough $n$'s}.
\end{equation}
\end{theorem}

We define the mass $J_n$ of the favourite exit point for the polymer,
$$
J_n= \max_{x \in \Z^d} \mu_{n-1}(\w_n=x)
$$
We view $J_n \in [0,1]$ as an index of localization of the polymer.
When $J_n$ vanishes, the polymer is delocalized in the sense that it spreads 
over all sites; This is the case for $\b=0$. 
On the other hand, when $J_n$ does not  vanish, the polymer is
strongly localized in the sense it has a significant probability to
go through a few special sites. More precisely,
\begin{definition} \label{def:loc}
The polymer is delocalized if
$$\lim_{n \to \8} J_n = 0 \qquad Q-{\rm a.s.}\;,$$
and localized if
$${\rm Cesaro-}\liminf_{n \to \8} J_n > 0  \qquad Q-{\rm a.s.}$$
where the $liminf$ is taken in the Cesaro sense, i.e. 
$${\rm Cesaro-}\liminf_{n \to \8} J_n=
\liminf_{n \to \8} \frac{1}{n} \sum_{t=1}^n J_t\;.$$
\end{definition}

In view of the relations 
$$J_n^2 \leq I_n \leq J_n\;,$$ and 
following \cite{CSY03},
it is not difficult to derive, 
 from Theorem
 \ref{th:equiv2} and Proposition \ref{prop:freeenergy},
the following corollary.
\begin{corollary} \label{cor:loc0}
We have the equivalences
$$
p=\lm \iff {\rm delocalization}
$$
and
$$
p<\lm \iff {\rm localization}
$$
\end{corollary}
In particular, for all $\b$, 
either delocalization occurs or localization occurs. In other words, there is 
another dychotomy: for all fixed $\b$, either $J_n$ vanishes for almost
every environment, or for almost
every environment, $\liminf J_n$ is positive in the Cesaro sense.
\medskip

Now, from Theorem \ref{th:WDa} and Proposition  \ref{th:WD}
we derive
\begin{corollary} \label{cor:loc}
(a) Assume  (\ref{hyp:2d}), and  either
(i) $d=1$ and $\a \in (0,1)$, or (ii) $d=2$, $\a \neq
2$, or (iii) $d\geq 3, \a \in (0,2]$. Then, delocalization holds 
for all $\b$ with (\ref{L2}).
\\
(b) Under the condition (\ref{eq:KP}),  localization holds. 
\end{corollary}

\small

\end{document}